\documentclass[12pt]{amsart}

\newcommand{\elle}[1]{L^{#1}(\Omega)}
\newcommand{\huz}{H^1_0(\Omega)}
\newcommand{\sob}[2]{W^{#1}_{#2}(\Omega)}

\newcommand{\enne}{\mathbb{N}}
\newcommand{\erre}{\mathbb{R}}

\newcommand{\io}{\int_{\Omega}}
\newcommand{\ik}{\int_{\{|\un| \geq k\}}}

\newcommand{\dive}{{\rm div}}
\newcommand{\norma}[2]{\|#1\|_{\lower 4pt \hbox{$\scriptstyle #2$}}}
\newcommand{\fn}{f_{n}}
\newcommand{\un}{u_{n}}

\newcommand{\vp}{\varphi}

\newcommand{\arrstre}{\renewcommand{\arraystretch}{2}}

\newcommand{\disp}{\displaystyle}
\newcommand{\ba}{\begin{array}}
\newcommand{\ea}{\end{array}}
\newcommand{\be}{\begin{equation}}
\newcommand{\ee}{\end{equation}}

\newcommand{\rife}[1]{(\ref{#1})}
\newtheorem{example}{\sc Example}

\newtheorem{ohss}[example]{\sc Remark}

\newtheorem{theo}[example]{\sc Theorem}

\begin{document}

\title{A semilinear problem with a $W^{1,1}_0$ solution}
\author{Lucio Boccardo, Gisella Croce, Luigi Orsina}
\address{L.B. -- Dipartimento di Matematica, ``Sapienza'' Universit\`{a} di Roma,
P.le A. Moro 2, 00185 Roma (ITALY)}
\email{boccardo@mat.uniroma1.it}
\address{G.C. -- Laboratoire de Math\'ematiques Appliqu\'ees du Havre, Universit\'e du Havre,
25, rue Philippe Lebon, 76063 Le Havre (FRANCE)}
\email{gisella.croce@univ-lehavre.fr}
\address{L.O. -- Dipartimento di Matematica, ``Sapienza'' Universit\`{a} di Roma,
P.le A. Moro 2, 00185 Roma (ITALY)}
\email{orsina@mat.uniroma1.it}



\begin{abstract}
We study a  degenerate elliptic equation, proving the existence of a  $W^{1,1}_0$ distributional solution.
\end{abstract}

\maketitle

In the study of elliptic problems, 
it is quite standard to find  solutions belonging either to $BV(\Omega)$ or to $W^{1,s}(\Omega)$, with $s>1$. 
 In this paper we prove the existence of a 
 $W^{1,1}_0$ distributional solution for the following
 boundary value problem:
\be\left\{
\arrstre
\ba{cl}
\disp
-\dive\bigg(\frac{a(x)\,\nabla u}{(1+b(x)|u|)^{2}}\bigg) + u = f & \mbox{in $\Omega $,}\\
\hfill u = 0 \hfill & \mbox{on $\partial\Omega $.}
\ea
\right.
\label{lineare}
\ee
Here   $\Omega$ is a bounded, open subset of $\erre^{N}$, with $N > 2$,
   $a(x),\,b(x)$ are  measurable functions such that
\begin{equation}
\label{ab}
0 <\alpha\leq a(x) \leq\beta,\quad  0 \leq b(x) \leq B,
\end{equation}
  with  $\alpha,\,\beta\in\erre^+$, $B \in \erre$ and
\begin{equation}
\label{f2}
f(x) \hbox{ belongs to } \elle2.
\end{equation}
 We are going to prove that problem \rife{lineare} has a distributional solution $u$ belonging to the non-reflexive Sobolev space $W^{1,1}_{0}(\Omega)$.

Problems like \rife{lineare} have been extensively studied in the past. In \cite{bdo}, existence and regularity results were obtained for
\be
\left\{
\arrstre
\ba{cl}
\disp
-\dive\bigg(\frac{a(x)\,\nabla u}{(1+|u|)^{\theta}}\bigg) = f & \mbox{in $\Omega $,}\\
\hfill u = 0 \hfill & \mbox{on $\partial\Omega$,}
\ea
\right.
\label{p-teta}
\ee
where $0<\theta\leq1$ and $f$ belongs to $\elle{m}$ for some $m \geq 1$.
A whole range of existence results was proved, yielding solutions belonging to some Sobolev space $\sob{1,q}0$, with $q=q(2,m)\leq 2$
or entropy solutions.
In the case where $\theta>1$ a non-existence result for constant sources has been proved in \cite{ABFOT}.

As pointed out in \cite{bb}, existence of solutions can be recovered for any value of $\theta>0$,
by adding a lower order term of order zero.
If we consider the problem
\be
\left\{
\arrstre
\ba{cl}
\disp
-\dive\bigg(\frac{a(x)\nabla u}{(1+|u|)^{2}}\bigg) + u = f& \mbox{in $\Omega $,}\\
\hfill u = 0 \hfill & \mbox{on $\partial\Omega$,}
\ea
\right.
\label{p-2}
\ee
with $f$ in $\elle{m}$, then  the following results can be proved  (see \cite{bb} and \cite{croce}):

\begin{enumerate}
\item[i)]
if $2 < m < 4$, then there exists a distributional solution in $ W^{1,\frac{2m}{m + 2}}_0(\Omega)\cap \elle{m}$;
\item[ii)]
if $1 \leq	m \leq 2$, then there exists an entropy solution in $\elle{m}$ whose gradient belongs to the Marcinkiewicz space $M^{\frac{m}{2}}(\Omega)$.
\end{enumerate}

In this paper we deal with the borderline case $m=2$, improving the above results as follows.

\begin{theo}\label{thm-lineare}\sl
Assume \rife{ab} and \rife{f2}. Then there exists a distributional
solution $u \in W^{1,1}_0(\Omega)\cap L^2(\Omega)$ to problem (\ref{lineare}), in the sense that
$$
\io \frac{a(x)\,\nabla u\cdot \nabla \varphi}{(1+b(x)|u|)^{2}} + \io u\,\varphi = \io f\,\varphi\,,
$$
for all $\varphi \in W^{1,\infty}_0(\Omega)$.
\end{theo}
\begin{ohss}\rm
If the operator is nonlinear with respect to the gradient,
existence of distributional solutions will be studied in a forthcoming paper (\cite{bco}).
\end{ohss}

\noindent{\sl Proof of Theorem 1.}

{\sl Step 1.}
We begin by approximating our boundary value problem \rife{lineare} and we
consider a sequence $\{\fn\}$ of $\elle\infty$ functions such that $\fn$ strongly converges to $f$ in $\elle2$, and $|\fn| \leq |f|$ for every $n$ in $\enne$. The same technique of \cite{bb} assures the existence of a solution $\un$ in $\huz \cap \elle\infty$ of
\begin{equation}\label{ppn_0}
\left\{
\arrstre
\ba{cl}
\disp
-\dive\bigg(\frac{a(x)\,\nabla \un}{(1+b(x)|\un|)^{2}}\bigg) + \un = \fn & \mbox{in $\Omega $,}\\
\hfill \un = 0 \hfill & \mbox{on $\partial\Omega $.}
\ea
\right.
\end{equation}
Indeed,
let $M_n = \norma{f_n}{\elle\infty} + 1$, and consider the problem
\be\label{pbbase}
\left\{
\arrstre
\ba{cl}
\disp
-\dive\bigg(\frac{a(x)\nabla w }{(1+b(x)|T_{M_n}(w)|)^{2}}\bigg) + w = f_n & \mbox{in $\Omega$,} \\
\hfill w = 0 \hfill & \mbox{on $\partial\Omega$,}
\ea
\right.
\ee
where $T_{k}(s) = \max(-k,\min(s,k))$ for $k \geq 0$ and $s$ in $\erre$. The existence of  a weak solution $w$ in $\huz$ of \rife{pbbase} follows from Schauder's theorem. Choosing $(|w| - \norma{f_n}{\elle\infty})_{+}\,{\rm sgn}(w)$ as a test function we obtain, dropping the nonnegative first term,
$$
\io |w|\,(|w| - \norma{f_n}{\elle\infty})_{+} \leq \io \norma{f_n}{\elle\infty}\,(|w| - \norma{f_n}{\elle\infty})_{+}\,.
$$
Thus,
$$
\io (|w| - \norma{f_n}{\elle\infty})\,(|w| - \norma{f_n}{\elle\infty})_{+} \leq 0\,,
$$
so that $|w| \leq \norma{f_n}{\elle\infty} < M_n$. Therefore, $T_{M_n}(w) = w$, and $w$ is a bounded weak solution of \rife{ppn_0}.

{\sl Step 2.} We prove some {\sl a priori} estimates on the sequence $\{u_n\}$.
Let $k \geq 0$, $i > 0$, and let $\psi_{i,k}(s)$ be the function defined by
$$
\psi_{i,k}(s) =
\left\{
\ba{cl}
0 & \mbox{if $0 \leq s \leq k$,}
\\
i(s-k) & \mbox{if $k < s \leq k + \frac1i$,}
\\
1& \mbox{if $s > k + \frac1i$,}
\\
\psi_{i,k}(s)=-\psi_{i,k}(-s) & \mbox{if $s < 0$.}
\ea
\right.
$$
Note that
$$
\lim_{i \to +\infty}\,\psi_{i,k}(s) = \left\{
\ba{cl}
1& \mbox{if $s > k$,}
\\
0 & \mbox{if $|s| \leq k$,}
\\
-1 & \mbox{if $s < -k$.}
\ea
\right.
$$
We choose $|\un|\,\psi_{i,k}(\un)$ as a test function in \rife{ppn_0}, and we obtain
$$
\arrstre
\ba{l}
\disp
\io \frac{a(x)|\nabla\un|^2}{(1 + b(x)|\un|)^{2}} |\psi_{i,k}(\un)|
+
\io \frac{a(x)|\nabla\un|^2}{(1 + b(x)|\un|)^{2}} \psi'_{i,k}(\un) |\un|
\\
\qquad
\disp
+
\io \un |\un| \psi_{i,k}(\un)
=
\io \fn |\un| \psi_{i,k}(\un)
\,.
\ea
$$
Since $\psi'_{i,k}(s) \geq 0$, we can drop the second term; using \rife{ab}, and the assumption $|\fn| \leq |f|$, we have
$$
\alpha\io \frac{|\nabla\un|^{2}}{(1 + b(x)|\un|)^{2}}\,|\psi_{i,k}(\un)|
+
\io \un |\un| \psi_{i,k}(\un)
\leq
\io |f| |\un| |\psi_{i,k}(\un)|
\,.
$$
Letting $i$ tend to infinity, we thus obtain, by Fatou's lemma (on the left hand side) and by Lebesgue's theorem (on the right hand side, recall that $\un$ belongs to $\elle\infty$),
\be\label{motherofall}
\alpha
\ik \frac{|\nabla\un|^{2}}{(1 + b(x)|\un|)^{2}}
+
\ik |\un|^{2}
\leq
\ik |f|\,|\un|
\,.
\ee
Dropping the nonnegative first term in \rife{motherofall} and using H\"older's inequality on the right hand side, we obtain
$$
\ik|\un|^{2}
\leq
\bigg[\ik |f|^{2}\bigg]^{\frac{1}{2}}
\,
\bigg[\ik |\un|^{2} \bigg]^{\frac{1}{2}}
\,.
$$
Simplifying equal terms we thus have
\be\label{aa}
\ik|\un|^{2} \leq \ik|f|^{2}\,.
\ee
For $k = 0$, \rife{aa} gives
\be\label{ll}
\io |\un|^{2} \leq \io |f|^{2}\,,
\ee
so that $\{\un\}$ is bounded in $\elle{2}$. This fact implies in particular that
\be\label{qq}
\lim_{k \to +\infty}\,{\rm meas}(\{|\un| \geq k\}) = 0\,,
\quad
\mbox{uniformly with respect to $n$.}
\ee

From \rife{motherofall}, written for $k = 0$, dropping the nonnegative second term and using that $b(x)\leq B$, we have
$$
\alpha \io\frac{|\nabla \un|^2}{(1 +B |\un|)^{2}}
\leq
\io |f|\,|\un|\,.
$$
H\"older's inequality on the right hand side then gives
$$
\alpha \io \frac{|\nabla \un|^2}{(1+B|\un|)^{2}}
\leq
\bigg[\io |f|^{2} \bigg]^{\frac{1}{2}}
\bigg[\io |\un|^{2} \bigg]^{\frac{1}{2}}\,,
$$
so that, by \rife{ll}, we infer that
\be\label{bb}
\alpha
\io\frac{|\nabla \un|^2}{(1+B| \un |)^{2}}
\leq
\io |f|^{2}\,.
\ee

{\sl Step 3.}
We prove that, up to subsequences, the sequence $\{u_n\}$ strongly converges in $\elle2$ to some function $u$.

From (\ref{bb}) we deduce that $v_{n} = \log(1 + B|\un|){\rm sgn}(\un)$ is bounded in $\huz$. Therefore, up to subsequences, it converges to some function $v$ weakly in $\huz$, strongly in $\elle2$, and almost everywhere in $\Omega$. If we define $u = \frac{{\rm e}^{|v|} - 1}{B}{\rm sgn}(v)$, then $\un$ converges almost everywhere to $u$ in $\Omega$. Let now $E$ be a measurable subset of $\Omega$; then
$$
\arrstre
\ba{r@{\hspace{2pt}}c@{\hspace{2pt}}l}
\disp
\int_{E}|\un|^{2}
& \leq &
\disp
\int_{E\cap\{|\un| \geq k\}}|\un|^{2}
+\int_{E\cap\{ |\un|<k\}}|\un|^{2}
\\
& \leq &
\disp
\ik|f|^{2}+k^{2}{\rm meas}(E)\,,
\ea
$$
where we have used \rife{aa} in the last passage. Thanks to \rife{qq}, we may choose $k$ large enough so that the first integral is small, uniformly with respect to $n$; once $k$ is chosen, we may choose the measure of $E$ small enough such that the second term is small. Thus, the sequence $\{\un^{2}\}$ is equiintegrable and so, by Vitali's theorem, $\un$ strongly converges to $u$ in $\elle2$.

{\sl Step 4.}
We prove that, up to subsequences, the sequence $\{u_n\}$ weakly converges to  $u$ in $\sob{1,1}0$.

Let again $E$ be a measurable subset of $\Omega$, and let $i$ be in $\{1,\ \ldots,\ N\}$. Then
$$
\arrstre
\begin{array}{rcl}
\disp
\int_{E} |\partial_{i}\un|
& \leq &
\disp 
\int_{E} |\nabla \un|
=
\int_{E} \frac{|\nabla \un|}{1 + B|\un|}\,(1 + B|\un|)
\\
&
\leq
&
\disp \bigg[
\int_{E}\frac{|\nabla \un|^2}{(1+B| \un |)^{2}}
\bigg]^\frac 12
\bigg[\int_{E} (1+B|\un|)^{2}\bigg]^\frac{1}{2}
\\
&
\leq
&
\disp
\bigg[\frac{1}{\alpha}\io |f|^2\bigg]^{\frac12}
\bigg[2{\rm{meas}}(E)+2 B^2\int_{E} |\un|^2\bigg]^{\frac12}\,,
\ea
$$
where we have used \rife{bb} in the last passage. Since the sequence $\{\un\}$ is compact in $\elle{2}$, we have that the sequence $\{\partial_{i}\un\}$ is equiintegrable. Thus, by Dunford-Pettis theorem, and up to subsequences, there exists $Y_{i}$ in $\elle1$ such that 
$\partial_{i} \un$ weakly converges to $Y_{i}$ in $\elle1$. Since $\partial_{i}\un$ is the distributional derivative of $\un$, we have, for every $n$ in $\enne$,
$$
\io \partial_{i} \un\,\vp = -\io \un\,\partial_{i} \vp\,,
\quad
\forall \vp \in C^{\infty}_{0}(\Omega)\,. 
$$
We now pass to the limit in the above identities, using that $\partial_{i}\un$ weakly converges to $Y_{i}$ in $\elle1$, and that $\un$ strongly converges to $u$ in $\elle2$; we obtain
$$
\io Y_{i}\,\vp = -\io u\,\partial_{i} \vp\,,
\quad
\forall \vp \in C^{\infty}_{0}(\Omega)\,, 
$$
which implies that $Y_{i} = \partial_{i} u$, and this result is true for every $i$. Since $Y_{i}$ belongs to $\elle1$ for every $i$, $u$ belongs to $W^{1,1}_{0}(\Omega)$, as desired. 

Note now that, since $s \mapsto \log(1+B s)$ is Lipschitz continuous on $\erre^{+}$, and $u$ belongs to $W^{1,1}_{0}(\Omega)$, by the chain rule we have
$$
\nabla [\log(1 + B|u|)\,{\rm sgn}(u)] = \frac{\nabla u}{1 + B|u|}\,,
\quad
\mbox{almost everywhere in $\Omega$.}
$$
Hence, from the weak convergence of $v_{n}$ to $v$ in $\huz$ we deduce that
\be\label{ok}
\lim_{n \to +\infty}\,\frac{\nabla \un}{1 + B|\un|} = \frac{\nabla u}{1 + B|u|}\,,
\quad
\mbox{weakly in $(\elle2)^{N}$.}
\ee

{\sl Step 5.} We now pass to the limit in the approximate problems \rife{ppn_0}.

Both the lower order term and the right hand side give no problems, due to the strong convergence of $\un$ to $u$, and of $\fn$ to $f$, in $\elle2$.

For the operator term we can write, if $\vp$ belongs to $W^{1,\infty}_0(\Omega)$,
\be\label{tt}
\io \frac{a(x)\,\nabla \un \cdot \nabla \vp}{(1 + b(x)|\un|)^{2}}
=
\io a(x)\,\frac{\nabla \un}{1 + B|\un|}\cdot \nabla \vp \,\frac{1 + B|\un|}{(1 + b(x)|\un|)^{2}}
\,.
\ee
In the last integral, the first term is fixed in $\elle\infty$, the second is weakly convergent in $(\elle2)^{N}$ by \rife{ok}, the third is fixed in $(\elle\infty)^{N}$, and the fourth is strongly convergent in $\elle2$, since is bounded from above by $1 + B|\un|$, which is compact in $\elle2$. Therefore, we can pass to the limit to have that
$$
\lim_{n \to +\infty}\,\io \frac{a(x)\,\nabla \un \cdot \nabla \vp}{(1 + b(x)|\un|)^{2}} = \io \frac{a(x)\,\nabla u \cdot \nabla \vp}{(1 + b(x)|u|)^{2}}\,,
$$
as desired.

\begin{ohss}\rm
Note that if $b(x) \geq b > 0$ in $\Omega$, then we can choose test functions $\vp$ in $\huz$. Indeed,
$$
0 \leq \frac{1 + B|\un|}{(1 + b(x)|\un|)^{2}} \leq \frac{1 + B|\un|}{(1 + b|\un|)^{2}} \leq C(B,b)\,,
$$
for some nonnegative constant $C(B,b)$, so that we can rewrite \rife{tt} as 
$$
\io \frac{a(x)\,\nabla \un \cdot \nabla \vp}{(1 + b(x)|\un|)^{2}}
=
\io a(x)\,\frac{\nabla \un}{1 + B|\un|}\cdot \frac{\nabla \vp \,(1 + B|\un|)}{(1 + b(x)|\un|)^{2}}
\,,
$$
with the first term fixed in $\elle\infty$, the second weakly convergent in $(\elle2)^{N}$, and the third strongly convergent in the same space by Lebesgue's theorem.
\end{ohss}


\end{document}